\documentclass[a4paper, portrait]{article}
\usepackage{amsmath,amsfonts,amssymb,amsthm,multirow}
\usepackage{graphicx}
\usepackage{grffile}
\usepackage{epstopdf,url}
\usepackage {epsfig}

\newtheorem{theorem}{Theorem}

\newtheorem {lemma}{Lemma}

\title{Optimality of constant inflow for a linear system with a bottleneck entrance}
\author{Guy Katriel\footnote{Email: katriel@braude.ac.il, Tel.: 972-4-9086443 }\\ Department of Mathematics, ORT Braude College,\\ Karmiel, Israel\\}

\date{}

\begin{document}

\maketitle

\begin{abstract}
Sadeghi et al. \cite{sadeghi1,sadeghi2} considered a bottleneck system with periodic inflow rate, and proved that a constant-rate input maximizes the time-averaged output rate among all periodic inflow rates. Here we provide a short and elementary proof of this result, without use of optimal control theory. The new approach developed here allows us to prove an extension of the result to the case of a general non-periodic inflow rate.

\end{abstract}

\section{Introduction}

A question of fundamental importance in the theory of process control is whether it is possible to improve the performance of a
system by operating it in a time-varying manner rather than with constant parameters. Indeed it has been established that in 
some circumstances periodic operation can lead to efficiency gains \cite{bayen,colonius,silveston}.

In recent studies Sadeghi et al. \cite{sadeghi1,sadeghi2} considered this question in the context of a simple model for a system with bottleneck entrance, which, as they note, can approximate the behavior of traffic systems, queues for security checks, and biological machines. 
In this model occupancy is described by the variable $x\in [0,1]$, and inflow rate is a given non-negative function $\sigma(t)$. Occupancy 
increases at a rate proportional to the inflow and the vacancy $1-x$, and output rate $w=\lambda x$ ($\lambda>0)$ is proportional to 
the occupancy, so that
\begin{equation}\label{de}x'(t)=\sigma(t)(1-x(t))-\lambda x(t).\end{equation}

When the inflow rate $\sigma$ is constant, the system converges to the steady state $x=\frac{\sigma}{\lambda+\sigma}$, so that the 
output rate converges to $w=\frac{\lambda\sigma}{\lambda+\sigma}$. If the inflow rate $\sigma(t)$ is a $T$-periodic function, $x(t)$ converges to the unique $T$-periodic solution $x_p(t)$ of (\ref{de}), so that the long-term time-averaged output rate is given by
\begin{equation}\label{w}w[\sigma]=\lambda \cdot \frac{1}{T}\int_0^T x_p(t)dt.\end{equation}
The question of interest is to choose $\sigma(t)$ so as to maximize this average output rate, for a given average input rate
\begin{equation}\label{sa}\bar{\sigma}=\frac{1}{T}\int_0^T \sigma(t)dt.\end{equation} 
In \cite{sadeghi1} this question is studied in the case in which $\sigma(t)$ describes switching between two values, and it is proved that in this case constant-rate input provides the optimal average output rate for a given average input rate. In other words, switching can never out-perform a constant input rate. In \cite{sadeghi2} it is proved that the same result holds more generally, that is for an arbitrary periodic function $\sigma(t)$. The proof is quite elaborate: it uses Pontryagin's maximum principle to reduce the problem to the case of switching case, and then a detailed study of the switching case is made.
 
In section \ref{periodic} of this letter we obtain a very short and elementary proof of the optimality of constant imput rate among all 
periodic input rates in the bottleneck system, using only elementary calculus, without the use of optimal control theory. In fact we show that the result follows 
from an identity, which thus quantifies the gap between the outputs obtained in for a general periodic input rate and that obtained for
the corresponding constant input rate with the same average.

 In section \ref{general} we use the new approach developed here to prove an extension of the result for the case in which $\sigma(t)$ is not periodic, showing that the 
long-term average output rate cannot exceed that obtained using a constant input rate with the same long-term average.

\section{Periodic input rate}
\label{periodic}

\begin{theorem}\label{main}
For any non-negative periodic $\sigma(t)$ with $\sigma|_{[0,T]} \in L^1[0,T]$, the time-averaged output rate (\ref{w}) satisfies 
$$w[\sigma]\leq w[\bar{\sigma}]=\frac{\lambda\bar{\sigma}}{\lambda+\bar{\sigma}},$$
with $\bar{\sigma}$ given by (\ref{sa}), and equality holds only when $\sigma(t)=\bar{\sigma}$ almost-everywhere.
\end{theorem}

Theorem \ref{main} follows immediately from the following identity:
\begin{lemma}\label{lem}
For any non-negative $\sigma \in L^1[0,T]$, if $x_p$ is the solution of (\ref{de}) satisfying $x_p(0)=x_p(T)$, we have
$$w[\sigma]=w[\bar{\sigma}]-\frac{1}{T}\int_0^T \left(x_p(t)-x^*\right)^2 (\lambda+\sigma(t))dt,$$
where $x^*=\frac{\bar{\sigma}}{\lambda+\bar{\sigma}}$.
\end{lemma}
Note that the integral on the right-hand side is non-negative, and vanishes only if $x_p(t)=\frac{\bar{\sigma}}{\lambda+\bar{\sigma}}$,
that is only if $\sigma(t)=\bar{\sigma}$. Therefore Theorem \ref{main} follows.

To prove Lemma \ref{lem}, we write (\ref{de}) in the form
\begin{equation}\label{de1}(\lambda+\sigma(t)) x(t)=\sigma(t)-x'(t)\end{equation}
and integrate both sides over $[0,\tau]$ ($\tau>0$), obtaining
\begin{equation}\label{i10}\int_0^{\tau} (\lambda+\sigma(t)) x(t)dt=\int_0^\tau \sigma(t)dt-(x(\tau)-x(0)).\end{equation}
We now multiply both sides of (\ref{de1}) by $x(t)$
$$(\lambda+\sigma(t)) x(t)^2=\sigma(t)x(t)-x(t) x'(t)$$
and integrate over $[0,\tau]$, to obtain
$$\int_0^\tau(\lambda+\sigma(t)) x(t)^2dt=\int_0^\tau \sigma(t)x(t)dt- \frac{1}{2}(x(\tau)^2-x(0)^2),$$
which together with (\ref{i10}) gives
\begin{equation}\label{i20}\int_0^\tau(\lambda+\sigma(t)) x(t)^2dt=\int_0^\tau \sigma(t)dt-\lambda \int_0^\tau x(t)dt-(x(\tau)-x(0))-\frac{1}{2}(x(\tau)^2-x(0)^2).\end{equation}
Assuming now that $x=x_p$ is the periodic solution of (\ref{de}), and taking $\tau=T$, (\ref{i10}),(\ref{i20}) become, after dividing 
both sides by $T$,
\begin{equation}\label{i1}\frac{1}{T}\int_0^{T} (\lambda+\sigma(t)) x_p(t)dt=\bar{\sigma},\end{equation}
\begin{equation}\label{i2}\frac{1}{T}\int_0^T(\lambda+\sigma(t)) x_p(t)^2dt=\bar{\sigma}-w[\sigma].\end{equation}
We now compute, using (\ref{i1}),(\ref{i2}),
$$\frac{1}{T}\int_0^T \left(x_p(t)-x^*\right)^2 (\lambda+\sigma(t))dt$$
$$= \frac{1}{T}\int_0^T(\lambda+\sigma(t))x_p(t)^2dt+{x^*}^2\cdot \frac{1}{T}\int_0^T(\lambda+\sigma(t))dt-2x^*\cdot \frac{1}{T}\int_0^T (\lambda+\sigma(t))x_p(t) dt$$
$$=\bar{\sigma}-w[\sigma]+{x^*}^2\cdot (\lambda+\bar{\sigma})-2x^* \bar{\sigma} $$
$$=\bar{\sigma}-w[\sigma]+ \frac{\bar{\sigma}^2}{(\lambda+\bar{\sigma})^2}\cdot (\lambda+\bar{\sigma})-\frac{2\bar{\sigma}}{\lambda+\bar{\sigma}}\cdot \bar{\sigma} $$
$$= \frac{\lambda\bar{\sigma}}{\lambda+\bar{\sigma}}-w[\sigma]=w[\bar{\sigma}]-w[\sigma],$$
and we have proved Lemma \ref{lem}.

\section{Extension to non-periodic input rate}
\label{general}

We now no longer assume that $\sigma(t)$ is periodic, but only that it is a locally-$L^1$ non-negative function.
We now define 
\begin{equation}\label{sa1}\bar{\sigma}=\limsup_{t\rightarrow \infty}\frac{1}{\tau}\int_0^\tau \sigma(t)dt.\end{equation}
We use the $\limsup$ because in general a limit will not exist. However for the important class of almost-periodic functions, the 
limit will exist \cite{levitan}. In particular, when $\sigma(t)$ is periodic, (\ref{sa1}) coincides with (\ref{sa}).

The long-time average output rate is at most 
\begin{equation}\label{wn}w[\sigma]=\lambda\cdot \limsup_{\tau\rightarrow\infty}\frac{1}{\tau}\int_0^\tau x(t),\end{equation}
where $x$ is any solution of (\ref{de}). The choice of solution does not matter, since for any two solutions
$x_1,x_2$, their difference will be a solution of the corresponding homogeneous equation, so that
$$\Big|\frac{1}{\tau}\int_0^{\tau }x_1(t)dt-\frac{1}{\tau}\int_0^\tau x_2(t)dt\Big|=|x_1(0)-x_2(0)|\frac{1}{\tau}\int_0^{\tau }e^{-\int_0^t (\lambda +\sigma (s))ds}dt$$
$$\leq |x_1(0)-x_2(0)| \frac{1}{\tau}\int_0^{\tau}e^{-\lambda t}dt\leq \frac{|x_1(0)-x_2(0)|}{\lambda\tau}\rightarrow 0
$$
as $\tau\rightarrow\infty$.
Note that (\ref{wn}) coincides with (\ref{w}) when $\sigma(t)$ is periodic.
We will prove
\begin{theorem}\label{main1}
	For any non-negative $\sigma \in L^1_{loc}[0,\infty]$ we have 
	\begin{equation}\label{thi}w[\sigma]\leq w[\bar{\sigma}]=\frac{\lambda\bar{\sigma}}{\lambda+\bar{\sigma}},\end{equation}
	with $\bar{\sigma}$ given by (\ref{sa1}).
\end{theorem}

To prove this theorem, we set $x^*=\frac{\bar{\sigma}}{\lambda+\bar{\sigma}}$, and compute, using (\ref{i10}),(\ref{i20}) 
$$0\leq \frac{1}{\tau}\int_0^\tau \left(x(t)-x^*\right)^2 (\lambda+\sigma(t))dt$$
$$= \frac{1}{\tau}\int_0^\tau(\lambda+\sigma(t))x(t)^2dt+{x^*}^2\cdot \frac{1}{\tau}\int_0^\tau(\lambda+\sigma(t))dt-2x^*\cdot \frac{1}{\tau}\int_0^\tau (\lambda+\sigma(t))x(t)dt$$
$$=\frac{1}{\tau}\int_0^\tau \sigma(t)dt-\lambda\cdot \frac{1}{\tau}\int_0^\tau x(t)dt-\frac{1}{\tau}(x(\tau)-x(0))-\frac{1}{2\tau}(x(\tau)^2-x(0)^2)$$
$$+{x^*}^2\cdot \frac{1}{\tau}\int_0^\tau(\lambda+\sigma(t))dt-2x^*\cdot\frac{1}{\tau}\left(\int_0^\tau \sigma(t)dt-(x(\tau)-x(0)) \right) $$
$$=(1-x^*)^2\cdot \frac{1}{\tau}\int_0^\tau \sigma(t)dt+\lambda {x^*}^2-\lambda\cdot \frac{1}{\tau}\int_0^\tau x(t)dt$$
$$+(2x^*-1)\cdot \frac{1}{\tau}(x(\tau)-x(0))-\frac{1}{2\tau}(x(\tau)^2-x(0)^2) $$
so that
$$\lambda\cdot \frac{1}{\tau}\int_0^\tau x(t)dt\leq (1-x^*)^2\cdot \frac{1}{\tau}\int_0^\tau \sigma(t)dt+\lambda {x^*}^2$$
$$+(2x^*-1)\cdot \frac{1}{\tau}(x(\tau)-x(0))-\frac{1}{2\tau}(x(\tau)^2-x(0)^2). $$
Taking the $\limsup$ as $\tau\rightarrow \infty$ on both sides, the last two terms vanish, since $x(\tau)\in [0,1]$, and we get
$$w[\sigma]\leq \bar{\sigma}(1-x^*)^2+\lambda {x^*}^2=\frac{\lambda\bar{\sigma}}{\lambda+\bar{\sigma}},$$
concluding the proof.

\end{document}